\documentclass[12pt]{article}

\newcommand{\p}{\partial}

\newtheorem{proposition}{Proposition}

\begin{document}
\begin{center}
{\bf\large A Covariant Poisson Deformation Quantization with Separation of
Variables up to the Third Order}\\

\bigskip

{\sc Alexander Karabegov}\\

{\it Abilene Christian University} \\

{\it E-mail:} alexander.karabegov@math.acu.edu
\end{center} 

\medskip

\begin{abstract} We give a simple formula for the operator $C_3$ of the
standard deformation quantization with separation of variables on a
K\"ahler manifold $M$. Unlike $C_1$ and $C_2$, this operator can not be
expressed in terms of the K\"ahler-Poisson tensor on $M$. We modify $C_3$
to obtain a covariant deformation quantization with separation of
variables up to the third order which is expressed in terms of the Poisson
tensor on $M$ and thus can be defined on an arbitrary complex manifold
endowed with a Poisson bivector field of type (1,1).
\end{abstract}

\bigskip

\section{Introduction}

\bigskip

Let $M$ be a Poisson manifold with the Poisson bivector field $\eta$. Then
$\{f,g\} = \langle \eta, df \wedge dg\rangle$ is a Poisson bracket on
$M$.  If $\eta$
is nondegenerate, its inverse $\omega$ is a symplectic form on $M$. Denote
by $C^\infty(M)[[\nu]]$ the space of formal series in $\nu$ with
coefficients from $C^\infty(M)$. As introduced in \cite{BFFLS}, a formal
differentiable deformation quantization on $M$ is an associative algebra
structure on $C^\infty(M)[[\nu]]$ with the $\nu$-linear and $\nu$-adically
continuous product $\ast$ (named star-product) given on $f,g\in
C^\infty(M)$ by the formula
\begin{equation} \label{E:star} 
f \ast g = \sum_{r = 0}^\infty \nu^r C_r(f,g), 
\end{equation} 
where $C_r,\, r\geq 0,$ are bidifferential operators on $M$, $C_0(f,g) =
fg$ and $C_1(f,g) - C_1(g,f) = i\{f,g\}$.  We adopt the convention that
the unit of a star-product is the unit constant. Two differentiable
star-products $\ast,\ast'$ on a Poisson manifold $(M,\{\cdot,\cdot\})$ are
called equivalent if there exists an isomorphism of algebras $B:
(C^\infty(M)[[\nu]],\ast) \to (C^\infty(M)[[\nu]],\ast')$ of the form $B=
1 + \nu B_1 + \nu^2 B_2 + \dots,$ where $B_r, r\geq 1,$ are differential
operators on $M$. The existence and classification problem for deformation
quantization was first solved in the non-degenerate (symplectic) case (see
\cite{DWL}, \cite{OMY}, \cite{F1} for existence proofs and \cite{F2},
\cite{NT}, \cite{D}, \cite{BCG}, \cite{X} for classification) and then
Kontsevich \cite{K} showed that every Poisson manifold admits a
deformation quantization and that the equivalence classes of deformation
quantizations can be parameterized by the formal deformations of the
Poisson structure.

If $M$ is a K\"ahler manifold, there exist special deformation
quantizations on $M$ such that the bidifferential operators $C_r$ defining
the star-product differentiate their first argument only in
antiholomorphic directions and the second argument only in holomorphic
ones. In \cite{CMP1} all the star-products with separation of variables on
an arbitrary K\"ahler manifold $M$ were completely described and
parameterized by the formal deformations of the original K\"ahler
structure on $M$. In what follows we will deal with the standard
deformation quantization with separation of variables which corresponds to
the trivial deformation of the K\"ahler structure. This standard
star-product was independently constructed in \cite{BW} with the use of
the properly tuned Fedosov's quantization scheme. Recently it was shown in
\cite{N} that all deformation quantizations with separation of variables
can be obtained via Fedosov's method. A star-product with separation of
variables on an arbitrary K\"ahler manifold was also obtained in \cite{RT}
by interpreting integral formulas of Berezin's quantization formally.

The coefficients of bidifferential operators $C_r$ of the standard
star-product with separation of variables $\ast$ on a K\"ahler manifold
$(M,\omega)$ in holomorphic local coordinates are polynomials in partial
derivatives of the K\"ahler metric tensor $g_{k\bar l}$ and its inverse
$g^{\bar lk}$. Notice that the tensor $g^{\bar lk}$ defines a global
Poisson bivector field of type (1,1) w.r.t. the complex structure on $M$.
Since the construction of the standard star-product with separation of
variables does not depend on the choice of local holomorphic coordinates,
the operators $C_r$ are "covariant" or "geometric". It follows almost
immediately from Fedosov's method that the operators $C_r$ can be
expressed in terms of the K\"ahler connection, its curvature and covariant
derivatives of the curvature. In \cite{E} Engli\v s calculated in a
covariant form the formal Berezin transform up to the third order of what
turns out to be the standard deformation quantization with separation of
variables. As it was noticed in \cite{RT}, unlike the operators $C_1$ and
$C_2$, the operator $C_3$ can not be expressed in terms of the Poisson
tensor $g^{\bar lk}$ only (there are always terms containing $g_{k\bar
l}$). Therefore, the formulas for the operators $C_r$ can not be used to
define a deformation quantization with separation of variables on an
arbitrary complex manifold endowed with a Poisson bivector field of type
(1,1) w.r.t. the complex structure.

In this letter we obtain a fairly simple covariant formula for the
operator $C_3$. Then we extract a non-covariant locally defined
one-differentiable (i.e., first order in both arguments) summand from
$C_3$ such that the rest is expressed in terms of $g^{\bar lk}$ only.
Finally, we modify $C_3$ by adding a covariant one-differentiable
bidifferential operator such that the sum is expressed in terms of
$g^{\bar lk}$ only. Thus we obtain a global covariant associative
star-product with separation of variables up to the third order on an
arbitrary complex manifold endowed with a Poisson bivector field of type
(1,1).

It is quite plausible that Kontsevich Formality can be established in the
"separation of variables" setting and a non-covariant Kontsevich-type star
product with separation of variables can be defined in terms of a
(possibly degenerate) Poisson tensor $g^{\bar lk}$ and then
globalized by the methods of \cite{CFT}. However, the construction of a
global star-product on a Poisson manifold from \cite{CFT} involves
non-canonical choices and is unlikely to produce directly a {\it
covariant} star-product with separation of variables. It would be
interesting to show the existence of such covariant star-products.

{\bf Acknowledgments} I am very grateful to M.~Engli\v s, B.~Fedosov, and
B.~Tsygan for interesting discussions. A part of this work had been done
when I was visiting the Mathematics department of the Pennsylvania State
University. I would like to thank them for their warm hospitality.

\section{The standard deformation quantization\\ 
with separation of variables}\label{S:sv}

A differentiable star-product $\ast$ on $M$ can be localized (restricted)
to any open subset $U\subset M$.  We shall retain the same notation $\ast$
for its restriction. For $f,g\in C^\infty(U)[[\nu]]$ denote by $L_f$ and
$R_g$ the operators of left star-multiplication by $f$ and of right
star-multiplication by $g$ respectively, so that $L_fg=f\ast g=R_gf$.
Denote by ${\cal L}(U)$ and ${\cal R}(U)$ the algebras of left and right
star-multiplication operators respectively. Notice that in order for a
product $\ast$ to be associative it is necessary and sufficient that for
any local functions $f,g$ the operators $L_f$ and $R_g$ commute,
$[L_f,R_g]=0$.

Let $M$ be a K\"ahler manifold with the K\"ahler form $\omega$. Its
inverse $\eta$ is a Poisson bivector field of type (1,1) w.r.t. the
complex structure on $M$. It determines a Poisson bracket on $M$. A formal
differentiable deformation quantization on $M$ is called quantization with
separation of variables if for any open subset $U \subset M$ and any
functions $a,b,f \in C^\infty(U)$ where $a$ is holomorphic and $b$
antiholomorphic, $a \ast f = af,\ f \ast b = bf$, that is, the operators
$L_a$ and $R_b$ are pointwise multiplication operators, $L_a = a, R_b =
b$. This means that the bidifferential operators $C_r$ defining the
star-product differentiate their first argument only in antiholomorphic
directions and the second argument only in holomorphic ones.

Recall the construction of the standard deformation quantization with
separation of variables on $M$ from \cite{CMP1}.  Let $(U,\{ z^k\})$ be an
arbitrary contractible holomorphic coordinate chart on $M$.  Locally
$\omega = - i g_{k\bar l} dz^k \wedge d\bar z^l$. Set $\p_k =
\frac{\p}{\p z^k}$ and $\bar\p_l = \frac{\p}{\p \bar z^l}$.  Pick a
K\"ahler potential $\Phi$ of the form $\omega$ on $U$, so that $g_{k\bar
l} = \p_k\bar\p_l\Phi$.  There exists a unique star-product with
separation of variables $\ast$ such that on any local chart $U$ as above
$L_{\p_k\Phi} = \p_k\Phi + \nu\p_k$ and $R_{\bar\p_l\Phi} = \bar\p_l\Phi +
\nu\bar\p_l$. The algebra of left star-multiplication operators ${\cal
L}(U)$ is described as the commutant of the operators $R_{\bar z^l} = \bar
z^l$ and $R_{\bar\p_l\Phi} = \bar\p_l\Phi + \nu\bar\p_l$ in the algebra of
all formal differential operators on $U$.  Once ${\cal L}(U)$ is known,
one can recover the product $f\ast g$ of any two functions $f,g \in
C^\infty(U)$ as follows. Since $f = f\ast 1 = L_f1$, then $L_f$ is found
as the unique operator $A\in {\cal L}(U)$ such that $A1 = f$. Now $f \ast
g = L_fg$. Recall how to find the operator $A = L_f$. Since $L_f$ commutes
with $R_{\bar z^l} = \bar z^l$, it does not contain antiholomorphic
partial derivatives. Since it commutes with $R_{\bar\p_l\Phi} =
\bar\p_l\Phi + \nu\bar\p_l$, one obtains recursive equations for the
components of $A = A_0 + \nu A_1 + \dots$,
\begin{equation}\label{E:rec}
[A_r,\bar\p_l\Phi] = [\bar\p_l, A_{r-1}],\quad r = 1,2,\dots 
\end{equation}
It was shown in \cite{CMP1} that equations (\ref{E:rec}) deliver a unique
operator $A$ commuting with antiholomorphic functions, such that $A1 = f$.
Finally we get that
\begin{equation}\label{E:ca}
C_r(f,g) = A_rg.
\end{equation}
The star-product $\ast$ is well defined globally on $M$. It does not
depend on the choices of local coordinates and K\"ahler potentials. The
coefficients of the bidifferential operators $C_r$ of the standard
star-product $\ast$ written in local holomorphic coordinates are
polynomials in partial derivatives of the K\"ahler metric tensor $g_{k\bar
l}$ and of its inverse $g^{\bar lk}$ (which determines the Poisson
bivector field $\eta = i g^{\bar lk}\p_k \wedge \bar \p_l$).
The operators $C_r$ are coordinate-independent and therefore "covariant"
or "geometric". They can be written in a "covariant" form in terms of the
K\"ahler connection on $M$.

\section{A covariant formula for the operator $C_3$}\label{S:c3}

First introduce and recall some notations and standard facts from K\"ahler
geometry. Throughout the paper we will use Einstein's summation
convention. Let $(U,\{ z^k\})$ be an arbitrary contractible holomorphic
coordinate chart on a K\"ahler manifold $(M,\omega)$ and $\Phi$ be a
K\"ahler potential of $\omega$ on $U$. For $A$ and $B$ a holomorphic and
antiholomorphic multi-indices respectively, set $g_{A\bar B} =
\p_A\bar\p_B\Phi$. In particular, $g_{k\bar l}=\p_k\bar\p_l\Phi$. The
Christoffel symbols of the K\"ahler connection $\nabla$ on $U$ are given
by the following formulas:

\begin{equation}\label{E:gamma}
\Gamma_{kp}^s = g_{kp\bar t}\, g^{\bar ts},\ 
\Gamma_{\bar l\bar q}^{\bar t} = g^{\bar ts}\, g_{s\bar q\bar l}. 
\end{equation}
Subsequently covariantly differentiating a function $f$ we obtain the
following symmetric tensors: $f_{/A} = \nabla_A f$ and $f_{/\bar B} =
\nabla_{\bar B}f$ (here $A$ and $B$ are a holomorphic and antiholomorphic
multi-indices respectively). In particular, $f_{/k}=\p_k f$ and $f_{/\bar
l}=\bar\p_l f$. The tensors $g_{k\bar l}$ and $g^{\bar lk}$ will be used
to lower and raise tensor indices. The Jacobi identity for the Poisson
tensor $g^{\bar lk}$ takes the form
\begin{equation}\label{E:jac}
g^{\bar lk}(\p_k g^{\bar nm}) = g^{\bar nk}(\p_k g^{\bar lm}),\
g^{\bar lk}(\bar\p_l g^{\bar nm}) = g^{\bar lm}(\bar\p_l g^{\bar nk}).
\end{equation}
The two following tensors obtained from the curvature tensor of the
K\"ahler connection $\nabla$
\[
 R_{p\bar q k\bar l} = g_{p\bar n}R^{\bar n}_{\bar qk\bar l} \mbox{\quad 
and \quad}
R_{k\bar l}^{\bar qp} = g^{\bar gm} R^p_{mk\bar l}
\]
are given by the formulas
\begin{equation}\label{E:curvdn}
R_{p\bar q k\bar l} 
= g^{\bar nm}\,g_{m\bar q\bar l}\,g_{pk\bar n} -
g_{p\bar q k\bar l}
\end{equation}
and
\begin{equation}\label{E:curvup}
R_{k\bar l}^{\bar qp} = \p_k \bar \p_l g^{\bar qp} - (\bar
\p_l g^{\bar qm})(\p_k g^{\bar np})g_{m\bar n}.
\end{equation}
Using recursive equations (\ref{E:rec}), formulas (\ref{E:jac}),
(\ref{E:curvdn}), (\ref{E:curvup}), and (\ref{E:ca}), one can derive the
following covariant formulas for the operators $C_1,C_2$ and $C_3$: 
\begin{equation}\label{E:c1}
C_1(\phi,\psi) = g^{\bar lk}(\bar\p_l\phi)(\p_k\psi) =
g^{\bar lk}\phi_{/\bar l}\,\psi_{/k},
\end{equation}
\begin{eqnarray}\label{E:c2}
C_2(\phi,\psi) = \frac{1}{2}(g^{\bar qp}g^{\bar
lk}(\bar\p_l\bar\p_q\phi)(\p_k\p_p\psi) + g^{\bar
qp}(\p_pg^{\bar lk})(\bar\p_l\bar\p_q\phi)(\p_k\psi) + \nonumber \\
(\bar\p_l g^{\bar qp})g^{\bar lk}(\bar\p_q\phi)(\p_k\p_p\psi) +
(\bar\p_l g^{\bar qp})(\p_p g^{\bar lk})(\bar\p_q\phi)(\p_k\psi)) =\\
\frac{1}{2}g^{\bar lk}g^{\bar qp}\phi_{/\bar
l\bar q}\,\psi_{/kp}.\nonumber
\end{eqnarray}
\begin{equation}\label{E:c3}
C_3(\phi,\psi) = \frac{1}{6}g^{\bar lk}g^{\bar qp}g^{\bar ts}\phi_{/\bar
l\bar q\bar t}\,\psi_{/kps} + \frac{1}{4} R^{\bar lk\bar qp}\phi_{/\bar
l\bar q}\psi_{/kp}.
\end{equation}

\section{A Poisson deformation quantization up to the third order}

If the formal parameter $\nu$ is nilpotent so that $\nu^{N+1} = 0$ for
some natural $N$, an associative product on a Poisson manifold $M$
determined by the formula
\[
f \ast g = \sum_{r = 0}^N \nu^r C_r(f,g) 
\]
with $C_0, C_1$ as in (\ref{E:star}) is called a "star-product up to the
$N$-th order".

It can be seen from formulas (\ref{E:c1}) and (\ref{E:c2}) in section
\ref{S:c3} that both operators $C_1$ and $C_2$ for the standard
deformation quantization with separation of variables can be expressed in
terms of the Poisson tensor $g^{\bar lk}$ only. However, this is not the
case for the operator $C_3$. We will call a bidifferential
operator regular if its coefficients in local holomorphic
coordinates can be written as polynomials in partial derivatives of the
tensor $g^{\bar lk}$. Introduce the following covariant global
bidifferential operators on $M$:
\begin{equation}\label{E:p}
P(\phi,\psi) := g^{\bar lk}g^{\bar qp}g^{\bar ts}\phi_{/\bar l\bar q\bar
t}\,\psi_{/kps}, \quad
Q(\phi,\psi) := - R^{\bar lk\bar qp}\phi_{/\bar l\bar q}\,\psi_{/kp},
\end{equation}
\[
\mbox{and} \quad R(\phi,\psi) := g^{\bar nm} R^{\bar lp}_{m\bar q}
\,R^{\bar qk}_{p\bar
n}\, \phi_{/\bar l}\,\psi_{/k}. 
\]
It turns out that all these operators coincide modulo regular operators.
Introduce the following locally defined (non-covariant) bidifferential
operator:
\begin{equation}\label{E:sing}
  S(\phi,\psi) := g_{m\bar n}(\bar \p_q g^{\bar
ls})(\p_s g^{\bar np})(\bar \p_t g^{\bar qm})(\p_p g^{\bar tk})
(\bar\p_l\phi)(\p_k\psi). 
\end{equation}

\begin{proposition}\label{P:prop}
All the operators $P,Q,R$ coincide with the operator $S$ modulo regular
operators.
\end{proposition}
{\it Proof.} We will prove only that the operator $R - S$ is regular. The 
rest of the proposition can be proved similarly. Using formula
(\ref{E:curvup}), rewrite $R(\phi,\psi)$ as follows:
\begin{eqnarray}\label{E:gup}
R(\phi,\psi) = 
g^{\bar nm}
(\p_m\bar\p_q g^{\bar lp} - (\bar\p_q g^{\bar la})
(\p_m g^{\bar bp})g_{a\bar b}) \nonumber\\
(\p_p\bar\p_n g^{\bar qk} - (\bar\p_n g^{\bar
qs})(\p_p g^{\bar tk})g_{s\bar t})(\bar\p_l \phi)(\p_k\psi).
\end{eqnarray}
Using formulas (\ref{E:jac}), we get
\begin{eqnarray}\label{E:first}
g^{\bar nm} 
(\bar\p_q g^{\bar la}) (\p_m g^{\bar bp}) g_{a\bar b}  =
(\bar\p_q g^{\bar la})g_{a\bar b}g^{\bar nm}(\p_m g^{\bar bp}) = \\
(\bar\p_q g^{\bar la})g_{a\bar b}g^{\bar bm}(\p_m g^{\bar np}) = 
(\bar\p_q g^{\bar lm})(\p_m g^{\bar np}). \nonumber
\end{eqnarray}
Similarly,
\begin{equation}\label{E:second}
g^{\bar nm} (\bar\p_n g^{\bar qs})(\p_p g^{\bar tk})g_{s\bar t} =
(\bar\p_n g^{\bar qm})(\p_p g^{\bar nk}).
\end{equation}
It follows from formulas (\ref{E:first}) and (\ref{E:second}) that
the operator $R(\phi,\psi)$ coincides modulo regular terms with the
operator
\begin{equation}\label{E:tilde}
\tilde S(\phi,\psi) = (\bar\p_q g^{\bar lm})(\p_m g^{\bar np}) (\bar\p_n
g^{\bar qs})(\p_p g^{\bar tk})g_{s\bar t}(\bar\p_l \phi)(\p_k\psi).
\end{equation}
We will now show that the operators $S$ and $\tilde S$ coincide.
Since $(\bar\p_n g^{\bar qs})g_{s\bar t} = - g^{\bar qs}g_{s\bar t\bar n}
= (\bar\p_t g^{\bar qs})g_{s\bar n},$ we get from formula (\ref{E:tilde})
that 
\begin{equation}\label{E:last}
\tilde S(\phi,\psi) = (\bar \p_q g^{\bar lm})(\p_m g^{\bar np})(\bar \p_t
g^{\bar qs})(\p_p g^{\bar tk}) g_{s\bar n} (\bar\p_l\phi)(\p_k\psi).
\end{equation}
It remains to swap the indices $s$ and $m$ in (\ref{E:last}) to obtain
formula (\ref{E:sing}) for the operator $S$.

In the notation of formulas (\ref{E:p})
formula (\ref{E:c3}) for the
operator $C_3$ can be rewritten as 
$$C_3 = \frac{1}{6} P - \frac{1}{4} Q.$$
It follows from Proposition \ref{P:prop} that the operator
\begin{eqnarray*}
\tilde C_3 (\phi,\psi)= C_3(\phi,\psi) + \frac{1}{12}R(\phi,\psi) =
\frac{1}{6}g^{\bar lk}g^{\bar qp}g^{\bar ts}\phi_{/\bar
l\bar q\bar t}\,\psi_{/kps} + \\
 \frac{1}{4} R^{\bar lk\bar qp}\phi_{/\bar
l\bar q}\psi_{/kp} + \frac{1}{12}
g^{\bar nm} R^{\bar lp}_{m\bar q} \,R^{\bar qk}_{p\bar
n}\, \phi_{/\bar l}\,\psi_{/k}
\end{eqnarray*}
is a regular covariant bidifferential operator that differs from the
operator $C_3$ by the one-differentiable (first order in both
arguments) operator $\frac{1}{12}R$. The following proposition immediately
follows
from the associativity of the standard star-product $\ast$ and
the fact that one-differentiable operators are Hochschild cocycles.
\begin{proposition} 
The formula
\[
  f\tilde\ast g = fg + \nu C_1(f,g) + \nu^2 C_2 (f,g) + \nu^3 \tilde
C_3(f,g)
\]
defines a regular covariant associative star-product with separation of
variables up to the third order on any K\"ahler manifold. It remains valid
on an arbitrary complex manifold endowed with a Poisson bivector field of
type (1,1) w.r.t. the complex structure.
\end{proposition}

\end{document}